\documentclass[11pt]{article}
\usepackage{amssymb,graphics,graphicx}
\newtheorem{definition}{{\bf Definition}}[subsection]

\newtheorem{theorem}{{\bf Theorem}}[subsection]
\newtheorem{proposition}{{\bf Proposition}}[subsection]

\newtheorem{lemma}{{\bf Lemma}}[subsection]
\newtheorem{remark}{{\bf Remark}}[subsection]

\let\ssection=\section
\renewcommand{\section}{\setcounter{equation}{0}\ssection}

\newcommand{\bbR}{\mathbb{R}}
\newcommand{\bbC}{\mathbb{C}}
\newcommand{\bbT}{\mathbb{T}}
\newcommand{\bbZ}{\mathbb{Z}}

\newcommand{\cH}{\mathcal{H}}
\newcommand{\cK}{\mathcal{K}}
\newcommand{\cI}{\mathcal{I}}
\newcommand{\cM}{\mathcal{M}}
\newcommand{\cO}{\mathcal{O}}
\newcommand{\cP}{\mathcal{P}}
\newcommand{\cQ}{\mathcal{Q}}
\newcommand{\cS}{\mathcal{S}}

\newcommand{\bA}{\mathbf{A}}
\newcommand{\bE}{\mathbf{E}}
\newcommand{\bS}{\mathbf{S}}
\newcommand{\bs}{\mathbf{s}}

\def\Ad{\mathop{\rm Ad}\nolimits}
\def\BT{\mathop{\bf BT}\nolimits}
\def\Coad{\mathop{\rm Coad}\nolimits}
\def\Conf{\mathop{\rm Conf}\nolimits}
\def\const{\mathop\mathrm{const.}\nolimits}

\def\Diff{\mathop{\mathrm{Diff}}\nolimits}
\def\Diffp{\mathop{\mathrm{Diff}_+}\nolimits}

\def\Div{\mathop{\rm Div}\nolimits}
\def\g{\mathop{\textrm{\rm g}}\nolimits}
\def\GF{\mathop{\bf GF}\nolimits}
\def\GLtwo{\mathop{\mathrm{GL}(2,\bbR)}\nolimits}
\def\Hess{\mathop{\mathrm{Hess}}\nolimits}
\def\id{\mathop{\mathrm{id}}\nolimits}
\def\Im{\mathop{\rm Im}\nolimits}

\def\Isom{\mathop{\rm Isom}\nolimits}
\def\PGL2{\mathop{\mathrm{PGL}(2,\bbR)}\nolimits}
\def\PSL2{\mathop{\mathrm{PSL}(2,\bbR)}\nolimits}
\def\tPSL2{\mathop{\mathrm{\widetilde{PSL}}(2,\bbR)}\nolimits}
\def\QS{\mathop{\mathrm{QS}}\nolimits}
\def\PSU11{\mathop{\mathrm{PSU}(1,1)}\nolimits}
\def\Rot{\mathop{\mathrm{Rot}}\nolimits}
\def\RP1{\mathop{\bbR P^1\!}\nolimits}
\def\RPn#1{\mathop{\bbR P^#1\!}\nolimits}
\def\S1{\bbT}
\def\SL2{\mathop{\mathrm{SL}(2,\bbR)}\nolimits}
\def\sl2{\mathop{\mathfrak{sl}(2,\bbR)}\nolimits}
\def\SO{\mathop{\mathrm{SO}}\nolimits}
\def\Vect{\mathop{\rm Vect}\nolimits}
\def\Vir{\mathop{\rm Vir}\nolimits}

\def\vfi{\varphi}
\def\vpi{\varpi}

\newcommand{\tS}{\widetilde{S}}
\newcommand{\halpha}{\widehat{\alpha}}

\def\x{\times}
\def\singleton#1{\lbrace#1\rbrace}
\def\class#1{\big[#1\big]}
\def\infinity{\singleton{\infty}}
\newcommand{\half}{\mathop{{\textstyle{\frac{1}{2}}}}\nolimits}
\newcommand{\threehalves}{\mathop{{\textstyle{\frac{3}{2}}}}\nolimits}
\newcommand{\QED}{\hskip 2truemm \vrule height3mm depth0mm width3mm}

\newcommand{\be}{\begin{equation}}
\newcommand{\ee}{\end{equation}}

\def\hfl#1#2{\smash{\mathop{\hbox to 12mm{\rightarrowfill}}
\limits^{\scriptstyle#1}_{\scriptstyle#2}}}

\def\vfl#1#2{\llap{$\scriptstyle#1$}\left\downarrow
\vbox to 6mm{}\right.\rlap{$\scriptstyle#2$}}

\def\diagram#1{\def\normalbaselines{\baselineskip=0pt
\lineskip=10pt\lineskiplimit=1pt}    \matrix{#1}}

\begin{document}

\baselineskip=13pt


\title{The Virasoro group \\and \\the fourth geometry of Poincar\'e}

\author{
{Christian DUVAL}
\thanks{E-mail: duval@cpt.univ-mrs.fr}\\
{\small    CPT-CNRS, Luminy Case 907}\\
{\small    F--13288 MARSEILLE Cedex 9 FRANCE}\\
            and\\
{Laurent GUIEU}
\thanks{E-mail: guieu@math.univ-montp2.fr}\\
{\small    GETODIM - D\'epartement de Math\'ematiques}\\
{\small    Universit\'e de Montpellier II}\\
{\small    F--34095 MONTPELLIER Cedex 5 FRANCE}
}   
\date{}

\maketitle

\begin{abstract}
We investigate, in some details, symplectic equivalence between several
conformal classes of Lorentz metrics on the hyper\-boloid of one
sheet $H^{1,1}\cong\bbT\x\bbT-\Delta$ and affine coadjoint orbits of the
group $\Diffp(\Delta)$ of orientation preserving diffeo\-morphisms of
$\Delta\cong\bbT$ with its natural projective structure.
This will allow for generalizations, namely, to the case of
arbitrary projective structures on null infinity.
\end{abstract}

\newpage


\section{Introduction}


According to the Riemann uniformization theorem, there exists
only three conformal types of simply connected Riemannian surfaces,
namely
$$
\begin{array}{c}
S^2\\[6pt]
K=1
\end{array}
\qquad
\qquad
\begin{array}{c}
\bbR^2\\[6pt]
K=0
\end{array}
\qquad
\qquad
\qquad
\begin{array}{c}
H^2\\[6pt]
K=-1.
\end{array}
$$

In the Lorentz case considered in this paper, the relevant geometry is the
so-called ``fourth'' geometry of Poincar\'e \cite{Poin} as opportunely
mentioned in \cite{KS2}, i.e., the Lorentz geometry of the hyperboloid of
one sheet
$$
\begin{array}{c}
H^{1,1}
\\[6pt]
K=\pm1.
\end{array}
$$

\begin{quote}
``\textsc{La quatri\`eme g\'eom\'etrie.} ---
Parmi ces axiomes implicites, il en est un qui semble m\'eriter
quelque attention, parce qu'en l'abandonnant, on peut construire une
quatri\`eme g\'eom\'etrie aussi coh\'erente que celle d'Euclide, de
Lobatchevsky et de Riemann.
[\dots]
Je ne citerai qu'un de ces th\'eor\`emes et je ne choisirai pas le plus
singulier~: \textit{une droite r\'eelle peut \^etre perpendiculaire \`a
elle-m\^eme}.''

\begin{flushright}
Henri Poincar\'e

La science et l'hypoth\`ese (1902)
\end{flushright}
\end{quote}

Let us, nevertheless, emphasize that a Lorentz uniformization
theorem is still not available, as of today---the problem lying
in the classification of the conformal boundaries \cite{Kul,Wei}.


This study has been triggered by previous work of Kostant and
Sternberg \cite{KS1,KS2} who first pointed out an intriguing relationship
between the Schwarzian derivative of a diffeomorphism of null infinity
$\bbT$ of the Lorentz hyperboloid $H^{1,1}$ and the transverse Hessian
of the conformal factor associated with this diffeomorphism (viewed as a
conformal transformation of~$H^{1,1}$). We contend that this
correspondence stems from a particular geometric object, namely the
cross-ratio as a four-point function associated with the canonical
projective structure of the projective line.

Such an observation prompted us to further investigate
the relationship between (i) the conformal geometry of the hyperboloid of
one sheet~$H^{1,1}$ and (ii) the Virasoro group, $\Vir$. 

Our contribution has therefore consisted in identifying several conformal
classes of Lorentz metrics on $H^{1,1}\cong\bbT^2-\Delta$ within the space
of projective structures on $\Delta\cong\bbT$, i.e., the (regular) dual of
$\Vect(\S1)$ \cite{Kiri}. In doing so, we have been able to give an
explicit, yet non standard, realization of the generic coadjoint orbits
\cite{Kiri,KY,Wit,Guieu1,Guieu2} of the Virasoro group in the framework of
$2$-dimensional real conformal geometry. Note that Iglesias
\cite{Igl} has also obtained other realizations of such orbits in quite a
different context.

\medskip
The paper is organized as follows.
\begin{itemize}

\item
Section \ref{TheLorentzHyperboloidOfOneSheet} describes in various ways
the Lorentz cylinder $\cH=\cS\x\cS-\Delta$ and its
associated conformal structure for special projective
structures of null infinity, i.e. the circle $\cS$. 

\item
In Section \ref{TheSchwarzianDerivative}, we briefly introduce the
Schwarzian $1$-cocycle $\bS$ of $\Diffp(\cS)$, while in Section
\ref{ConformalTransformations}, we recall the Kostant-Sternberg Theorem
\cite{KS2} and the basic notions attached to conformal Lorentz
structures on surfaces. 

\item
Our main results are presented in Section
\ref{SymplecticStructuresOnConformalClassesOfMetrics} where special,
infinite-dimensional, conformal classes of metrics $\g$ on $\cH$
are shown to be symplecto\-morphic to coadjoint orbits of the
group~$\Vir$---central extension of $\Conf_+(\cH)\cong\Diff_+(\cS)$. The
$\Conf_+(\cH)$-orbit of the flat Lorentz metric on the cylinder
corresponds to a zero central charge orbit, whereas the central charge $c$
of the other generic $\Vir$-orbits we investigate is related to the
(constant) curvature $K$ of $(\cH,\g)$ by $cK=1$. We, likewise, derive
the Bott-Thurston cocycle within the same framework.

\item
Some perspectives are finally drawn in Section \ref{ConclusionAndOutlook}.
It is, in particular, expected that our results allow for generalizations
that would, e.g., relate Kulkarni's Lorentz surfaces and
universal Teichm\"uller space.

\end{itemize}

\bigskip\noindent
\textbf{Acknowledgments:}
It is a pleasure for us to acknowledge enlightening conversations with
V.~Ovsienko and P.~Iglesias during the preparation of this article. We
would also like to thank H.~Heyer and J.~Marion, for the nice organization
of the Colloquium on \textit{Analysis on Infinite-Dimensional Lie Groups
and Algebras} held at CIRM in September 1997; a great many thanks for
their unwavering patience. 

\goodbreak

\section{The Lorentz hyperboloid of one sheet}
\label{TheLorentzHyperboloidOfOneSheet}

\subsection{An adjoint orbit in $\sl2$}

The single sheeted hyperboloid $H^{1,1}_c\hookrightarrow\bbR^{2,1}$ defined
for
$c\in\bbR^*_+$ by
\be
x^2+y^2-t^2=c
\label{HypEq}
\ee
carries a canonical Lorentz metric\footnote{In the physics literature
$H^{1,1}_1$ is called anti-de Sitter spacetime.} given by the
induced quadratic form 
\be
\g_c=dx^2+dy^2-dt^2.
\label{inducedMetric}
\ee
\begin{proposition}[\cite{Kul}, \cite{Wolf}]\label{Killing}
The hyperboloid of one sheet
$
H^{1,1}_c\cong\bbR\x\bbT
$
with radius $r=\sqrt{c}\neq0$ is the homogeneous space
$$
H^{1,1}_c=\SL2/\SO(1,1)
$$
which is symplectomorphic to the $\SL2$-adjoint orbit of
$$
\pmatrix{r&0\cr0&-r}\in\sl2.
$$
As a Lorentz manifold, $H^{1,1}_c$ is a space form of constant
curvature\footnote{Since $\g\longrightarrow-\g$ yields $K\longrightarrow-K$
and preserves the Lorentz signature $(+,-)$, we will admit $c<0$
in (\ref{curvature}); see \cite{Ghys2}. Recall that $K=\frac{1}{2}R$ where
$R$ is the scalar curvature.} 
\be
K=\frac{1}{c}
\label{curvature}
\ee
whose group of direct isometries is $\PSL2$.
\end{proposition}

\medskip
\begin{remark}\label{spaceOfOrientedGeodesicsOfH2}
{\rm 
The unit hyperboloid $H_1^{1,1}$ is also symplectomorphic to the manifold
of oriented geodesics of the Poincar\'e disk $H^2\cong\SL2/\SO(2)$. 
}
\end{remark}

From now on we will write $H$ as a shorthand notation for $H^{1,1}$
provided no confusion occurs.
\bigskip

\goodbreak

The following expression for the Lorentz metric (\ref{inducedMetric}) on
$H$ will prove useful. In view of (\ref{HypEq}), write 
$x = \varrho\sin\theta$, 
$y = \varrho\cos\theta$, 
$r = \varrho\sin\phi$, 
$t = \varrho\cos\phi$ 
so that the metric (\ref{inducedMetric}) takes the form
$\g_c=r^2\csc^2\!\phi\,(d\theta^2-d\phi^2)$. Putting now
$\theta_1=\theta+\phi$ and $\theta_2=\theta-\phi$, we obtain 
\be
\g_c = 
\frac{4c\,d\theta_1d\theta_2}{\left|e^{i\theta_1}-e^{i\theta_2}\right|^2}
\label{KostantSternberg1}
\ee
with (see (\ref{curvature}))
\be
c\in\bbR^*,
\label{cneq0}
\ee
yielding the canonical Killing metric on the
hyperboloid
\be
H\cong{}\bbT\x\bbT-\Delta
\label{KostantSternberg2}
\ee
globally parametrized by
$\theta_1,\theta_2\in\bbT=\bbR/(2\pi\bbZ)$ with
$\theta_1\neq\theta_2$. See, e.g.,~\cite{KS1}.
The transverse null foliations $\theta_1=\const$ and
$\theta_2=\const$ correspond to the rulings of the hyperboloid,
and the diagonal $\Delta$ is the conformal boundary \cite{Kul} (or
null infinity \cite{Penr}) of~$H$.

\subsection{
The Cayley-Klein model}\label{KleinSection}

The material of this Section has been borrowed from \cite{Cartan3} with
a slight adaptation to our framework.
\begin{definition}
{\rm 
An \textit{involution} of $\RP1$ is an homography $s\in\PGL2$ such that
$s^2=\id$ and $s\neq\id$. We will denote $\cI$ the space of involutions.
}
\end{definition}
In the projective plane $P$ associated to the vector space $\sl2$, there is
a distinguished conic $C$, defined by the light cone.

\begin{lemma}
The space of involutions is naturally identified with $P-C$.
\end{lemma}
The determinant map $\det:\GLtwo\longrightarrow\bbR^*$ descends, after
projectivization, as a map $\delta:\PGL2\longrightarrow\bbZ/(2\bbZ)$, that
defines the two connected components of the projective group. Then, we can
define $\cI_+ =\cI\cap \delta^{-1}(1)$ the space of direct involutions and
$\cI_- =\cI\cap \delta^{-1}(-1)$  the space of anti-involutions.
Let us denote by $D$ the interior of the convex hull of $C$ and by $\cK$
the complement of $D\cup C$ in $P$.

\begin{proposition} 
The space of direct involutions is naturally isomorphic to the disk $D$
and the space of anti-involutions to $\cK$.
\end{proposition}

\begin{remark}
\rm{
Topologically, $\cK$ is a M\"obius band.
}
\end{remark}

\begin{proposition}
The $2$-fold covering of orientations for $\cI_-$ is $C\times{}C-\Delta$.
The restriction of the projection
$\pi:\sl2\!-\singleton{0}\longrightarrow{}P$ to the Lorentz hyperboloid
$H$ is a $2$-fold covering on $\cK$.
\end{proposition}

There exists an isomorphism $P\cong\RPn{2}$ such that the conic
$C$ is mapped onto the unit circle $\S1$ in the affine plane
$\singleton{t=1}$, where $x,y,t$ are homogeneous coordinates in $\bbR^3$.
This isomorphism is given by the map
$$
X = 
\pmatrix{a&b\cr{}c&-a}
\longmapsto\frac{1}{2}
\pmatrix{2a\cr{}b+c\cr{}b-c}.
$$
Thus, we verify that the light
cone, whose equation is given by 
$\det(X)=0$, is mapped onto the conic of homogeneous equation
$x^2 + y^2 - t^2 = 0$.

  
In the Klein model, the complement
$
\cK=\{z\in\bbC\,|\,|z|>1\}
$
of the closed unit disk thus represents the projectivized hyperboloid
$P(H)$ in $\bbR P^2\cong{}P(\sl2)$. It is the space of geodesics of the
open unit disk, i.e., of the hyperbolic plane in the Klein model. See
Remark~\ref{spaceOfOrientedGeodesicsOfH2}.
\subsection{Projective structures}

In order to gain some insight into the preceding results, let us briefly
recall the notion of projective structure \cite{Cartan1,Cartan2,CF,Wil}. To
that end, we need the
\begin{definition}\label{projStructDef}
{\rm
A \textit{projective structure} $\vpi$ on a $n$-dimensional connected
manifold~$\cM$ is given by the following data:
\begin{enumerate}
\item
an immersion
$\Phi:\widetilde{\cM}\longrightarrow\bbR{}P^n$ defined on the universal
covering $\widetilde{\cM}$ of~$\cM$,
\item
a homomorphism 
$T:\pi_1(\cM)\longrightarrow\mathrm{PSL}(n+1,\bbR)$ 
\end{enumerate}
such that
\be
 \forall a \in\pi_1(\cM)
 \quad
 \Phi\circ a=T(a)\circ\Phi.
 \label{holo}
\ee
One calls $\Phi$ the \textit{developing map} and $T$ the 
\textit{holonomy} of the structure.
}
\end{definition}

We denote by $\vpi = [\Phi,T]$ the associated projective structure.
The developing map and the holonomy characterizes the structure up to
conjugation by the projective group, that is:
$$
\forall A \in \mathrm{PGL}(n+1,\bbR)
\quad
[A\circ\Phi,A\cdot{}T\cdot{}A^{-1}] = [\Phi,T].
$$

Such a structure is equivalently given by an atlas of projective charts 
$\vfi_i: U_i \subset \cM \longrightarrow \bbR{}P^n$ with transition
diffeomorphisms in $\mathrm{PGL}(n+1,\bbR)$.

In the $1$-dimensional case under study, and, more particularly in the case
of the circle~$\cS$, a projective structure $\vpi$ is given by a pair
$(\Phi,M)$ with $\Phi:\bbR\longrightarrow\RP1$ an
immersion and $M\in\PSL2$. Condition $\ref{holo}$
then reads
$$
\Phi(\theta+2\pi) = M\cdot \Phi(\theta).
$$

\goodbreak

It is a classic result \cite{Segal2,Guieu2} that the space $\cP(\cS)$ of
all projective structures on~$\cS$ is an affine space modeled on the space
$\cQ(\cS)$ of quadratic differentials $q=u(\theta)\,d\theta^2$ of~$\cS$.
The projective atlas associated with $q$ is obtained by locally solving the
third order non-linear differential equation
$
q=S(\Phi)
$
where $S$ stands for the Schwarzian derivative (see below).

From now on, we restrict considerations
to either choices of projective structures on $\cS$, namely
\begin{enumerate}
\item
the torus $\bbT=\bbR/(2\pi\bbZ)$ defined by the following developing
map\footnote{We use the notation $[z]=\bbR{}z$ for all
$z\in\bbC-\singleton{0}$.} (with trivial holonomy)
\be
\Phi(\theta)=[e^{i\theta}]
\qquad
\mbox{or}
\qquad
\Phi(\theta)=2\tan\frac{\theta}{2},
\label{T}
\ee

\item
the projective line $\RP1$ defined by the developing map
\be
\Phi(\theta)=\tan\theta
\qquad
\mbox{or}
\qquad
\Phi(t)= t.
\label{P1}
\ee
\end{enumerate}

\subsection{Lorentzian metric and cross-ratio}

Let us describe, following Ghys \cite{Ghys2},  how the canonical Lorentz
metric (\ref{KostantSternberg1}) on anti-de Sitter space 
(\ref{KostantSternberg2}) indeed originates from the cross-ratio
\be
(z_1,z_2,z_3,z_4) = \frac{(z_1-z_3)(z_2-z_4)}{(z_1-z_4)(z_2-z_3)}
\label{crossRatio}
\ee
of four points on the projective line \cite{Cartan2}.

Let us fix $(\theta_1,\theta_2)\in\bbT^2-\Delta$ and consider then a nearby
point $(\theta_3,\theta_4)=(\theta_1+d\theta_1,\theta_2+d\theta_2)$.
Put $z_j=e^{i\theta_j}$ for $j=1,\ldots,4$ and perform a Taylor
expansion of the cross-ratio (\ref{crossRatio}) at
$(\theta_1,\theta_2)$, so that
\goodbreak
\begin{eqnarray*}
(z_1,z_2,z_3,z_4) 
&=& 
\frac{
e^{i\theta_1}e^{i\theta_2}
\left(1-e^{id\theta_1}\right)\left(1-e^{id\theta_1}\right)
}
{
\left(e^{i\theta_1}-e^{i(\theta_2+d\theta_2)}\right)
\left(e^{i\theta_2}-e^{i(\theta_1+d\theta_1)}\right)
}\\
\\
&=&
\frac{
\left(-id\theta_1\right)\left(-id\theta_2\right)
}
{
\left(e^{i\theta_1}-e^{i\theta_2}\right)
\left(e^{i\theta_2}-e^{i\theta_1}\right)
e^{-i\theta_1}e^{-i\theta_2}
}
+\cdots\\
\\
&=&
\frac{
-d\theta_1 d\theta_2
}
{
\left|e^{i\theta_1}-e^{i\theta_2}\right|^2
}
+\cdots
\end{eqnarray*}
where the ellipsis ``$\cdots$'' stands for ``terms of order $\geq3$''.
One can thus claim that, up to higher order terms, the metric 
(\ref{KostantSternberg1}) on the unit hyperboloid $H$
(\ref{KostantSternberg2}) is given by
$
\g_1
=
-4\left(z_1,z_2,z_1+dz_1,z_2+dz_2\right)+\cdots
$ or, equivalently, by
\be
\g_1
=
-4\,\lim_{\varepsilon\rightarrow0}{
\frac{1}{\varepsilon^2}
\left(z_1,z_2,z_1+\varepsilon{}dz_1,z_2+\varepsilon{}dz_2\right)
}
\label{theMetric}
\ee
which is therefore conspicuously $\PSL2$-invariant.


\begin{figure}[h]
\begin{center}
\includegraphics[scale=0.5]{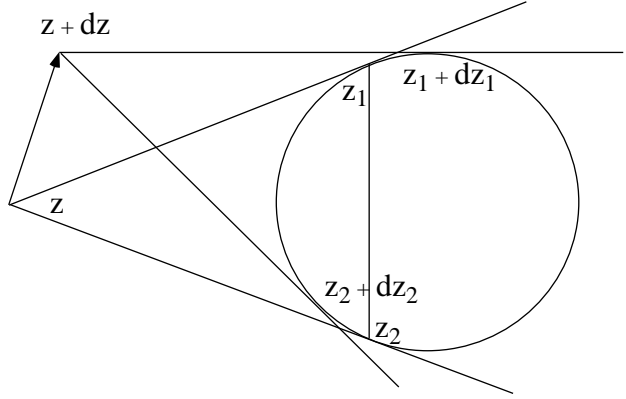}
\caption{\label{Klein}\textit{The Klein model}}
\end{center}
\end{figure}


\medskip

Resorting to Definition \ref{projStructDef}, we then have the 
\begin{theorem}
Consider the hyperboloid $\cH=\cS\times\cS-\Delta$ where the circle~$\cS$
has a projective structure defined by
$\Phi\in\Diff_\mathrm{loc}(\bbR,\RP1)$ as in (\ref{T}) or (\ref{P1}).
Then,
$\cH$ carries a natural $\PSL2$-invariant metric of the form
\be
\g_1 = 
(\Phi\times\Phi)^*\frac{4\,dt_1dt_2}{(t_1-t_2)^2}.
\label{generalg1}
\ee
\end{theorem}
\textit{Proof:}
The cross-ratio (\ref{crossRatio}) is $\PSL2$-invariant and so is the
Lorentz metric $4dt_1dt_2/(t_1-t_2)^2$ of $\RP1\x\RP1-\Delta$ given by
(\ref{theMetric}) with $z_j=t_j$ (see (\ref{P1})). In any cases (\ref{T})
or (\ref{P1}), the metric (\ref{generalg1}) defined on $\bbR^2-\Gamma$
where $\Gamma=(\Phi\times\Phi)^{-1}(\Delta)$ is automatically
$\pi_1(\cS)$-invariant thanks to (\ref{holo}). It is invariant, as well,
under the universal covering $\tPSL2$ of $\PSL2$. Hence, this metric
descends to $\cH=\cS\times\cS-\Delta=\pi\x\pi(\bbR^2-\Gamma)$ where
$\pi:\bbR\to\cS$ is the universal covering map. The projected metric
$\g_1$ is then clearly $\PSL2$-invariant. 
\QED

\medskip
Example (\ref{KostantSternberg1})
corresponds to the developing maps (\ref{T}); as for the first developing
map in (\ref{P1}), it leads via (\ref{generalg1}) to the metric of the
Klein model of Section \ref{KleinSection} (see Figure \ref{Klein}).

\section{The Schwarzian derivative}\label{TheSchwarzianDerivative}

\subsection{Osculating homography of a diffeomorphism}

Let $\vfi:\RP1\longrightarrow\RP1$ be a diffeomorphism and let $t_0\in\RP1$. We want
to find the homography $h\in\PGL2$ that best approximates the
diffeo\-morphism~$\vfi$ at this point $t_0$.
\begin{proposition}
This homography $h$ exists and is unique. It is completely defined by the
conditions
\begin{eqnarray*}
\hfil h(t_0)    &=& \vfi(t_0),\\
\hfil h'(t_0)   &=& \vfi'(t_0),\\
\hfil h''(t_0)  &=& \vfi''(t_0).
\end{eqnarray*}
\end{proposition}
The diffeomorphism $h^{-1}\circ\vfi$ has the $2$-jet of
the identity at $t_0$. The difference between $h$ and $\vfi$
starts, hence, at the third order derivative. (See, e.g.,~\cite{Ghys2}.)
\begin{definition}
{\rm 
The \textit{Schwarzian derivative} of $\vfi$ at the point $t_0$ is
$$
S(\vfi)(t_0) := \left( h^{-1} \circ \vfi \right) '''(t_0).
$$
}
\end{definition}
The quantity $S(\vfi)(t_0)$ measures how much does the diffeomorphism
$\vfi$ differ from an homography at the point $t_0$. All projective
information about $\vfi$ is encoded into the Schwarzian derivative. If we
identify the real projective line with $\bbR\cup\infinity$ by:
$[x,y] \longmapsto t=y/x$, we obtain the classical formula:
\be
S(\vfi) =
\left(
{\vfi'''(t)\over\vfi'(t)}-{3\over2}{\vfi''(t)^2\over\vfi'(t)^2}
\right) dt^2.
\label{TheFormule}
\ee
The graph $\Gamma_\vfi$ of our diffeomorphism is a simple closed curve
on $\RP1\times\RP1$.
\begin{definition}
{\rm 
The homography $h$ and its graph $\Gamma_h$ are respectively called the
\textit{osculating homography} and the \textit{osculating hyperbola} of
$\vfi$ at $t_0$.
}
\end{definition}

\subsection{The Schwarzian as a projective differential invariant}

\begin{theorem}[\cite{Green}]
The Schwarzian derivative is a third-order complete differential invariant for
the group of diffeomorphisms of the projective line.
\end{theorem}
More precisely, if $\vfi$ and $\psi$ are two diffeomorphisms of $\RP1$,
then
$$
S(\vfi) = S(\psi)
\qquad
\Leftrightarrow
\qquad
\exists A \in \PSL2,\ \psi = A \circ \vfi.
$$
\begin{theorem}[\cite{Bott,Kiri,Roger,Segal}]
The Schwarzian $S$ given by (\ref{TheFormule}) is a non trivial
$1$-cocycle, i.e.,
$$
S(\vfi\circ\psi) = \psi^*S(\vfi) + S(\psi)
\qquad
\forall\vfi,\psi\in\Diff_+(\RP1),
$$
on the group of orientation-preserving diffeomorphisms of $\RP1$ with
values in the $\Diff_+(\RP1)$-module of real quadratic differentials
$\cQ(\RP1)$ of of $\RP1$. Its kernel is $\PSL2$.
\end{theorem}

\begin{remark}
{\rm The Schwarzian cocycle (\ref{TheFormule}) is uniquely characterized
(up to a constant factor) by the property of having kernel $\PSL2$. 
}
\end{remark}

\subsection{Cartan formula of the cross-ratio}
\label{CartanFormulaSection}

A useful means for calculating the
Schwarzian derivative of a smooth map of the projective line is given by

\begin{theorem}[\cite{Cartan2}]\label{CartanThm}
Consider a smooth map $\vfi:\bbR P^1\longrightarrow\bbR P^1$ and
four points
$t_1,\ldots,t_4\in\bbR P^1$ tending to $t\in\bbR P^1$;
putting
$\tau_j=\vfi(t_j)$ one has
\be
\frac{(\tau_1,\tau_2,\tau_3,\tau_4)}{(t_1,t_2,t_3,t_4)}-1
=
\frac{1}{6}
S(\vfi)(t) (t_1-t_2)(t_3-t_4)+\hbox{\rm [higher order terms]}
\label{SchwarzCartan}
\ee 
where $S(\vfi)$ denotes the Schwarzian derivative (\ref{TheFormule}) of
$\vfi$.
\end{theorem}

This expression still makes sense for any smooth map of the
circle $\cS$ endowed with some projective structure given, for example, by
(\ref{T}) or (\ref{P1}). We, indeed, have the

\begin{definition}
{\rm
Let $\vfi:\cS\longrightarrow\cS$ be a smooth map identified with one
of its representatives\footnote{Choose any element of
$C^\infty(\bbR)$ that commutes with $\pi_1(\cS)$.} in
$C^\infty_{\pi_1(\cS)}(\bbR)$, then the \textit{Schwarzian} of $\vfi$ is
the pull-back of the Schwarzian (\ref{SchwarzCartan}) of the induced map
$\tilde\vfi$ of
$\RP1$, namely
\be
\bS(\vfi)=\Phi^*S(\tilde\vfi).
\label{generalSchwarzian}
\ee
}
\end{definition}

We note that (\ref{generalSchwarzian}) yields a well-defined quadratic
differential on $\cS$ since one trivially finds $a^*\bS(\vfi)=\bS(\vfi)$ in
view of $T(a)^*S(\tilde\vfi)=S(\tilde\vfi)$ for all
$a\in\pi_1(\cS)\cong\bbZ$.

\begin{proposition}
One has, locally,
\be
\bS(\vfi)=S(\vfi)+\vfi^*S(\Phi)-S(\Phi).
\label{generalSchwarzianLoc}
\ee
\end{proposition}
\textit{Proof:} 
Using $\tilde\vfi\circ\Phi=\Phi\circ\vfi$, one easily finds
$\Phi^*S(\tilde\vfi)(\theta)
=
S(\vfi)(\theta)+S(\Phi)(\vfi(\theta))\,\vfi'(\theta)^2
-S(\Phi)(\theta)$. 
\QED

\goodbreak

\section{Conformal transformations}\label{ConformalTransformations}

\subsection{Conformal Lorentz structures}

Let us recall some basic definitions and facts about $2$-dimensional
Lorentzian conformal geometry.
\begin{definition}[\cite{Kul}]\label{DefConf}
{\rm 
A \textit{conformal Lorentz structure} on a surface~$\Sigma$ is
characterized by a pair of transverse foliations; in other words, it is
given by a splitting
\be
T\Sigma=T_1\Sigma\oplus{}T_2\Sigma
\label{splitting}
\ee
into two trivial line bundles (light-cone field).
We call $N_1$ and $N_2$, respective\-ly, the spaces of leaves of the
two foliations of $\Sigma$.
}
\end{definition}

The leaves composing the ``grid'' associated to these foliations are,
locally, given by
$$
N_1:\theta_1=\const, \qquad N_2:\theta_2=\const
$$
The conformal structure is characterized by the global inter\-section
properties of the (null) leaves of $N_1$ and $N_2$.

One can associate to the splitting (\ref{splitting}) a class of
metrics on $\Sigma$, locally, of the form
$\g=F(\theta_1,\theta_2)\,d\theta_1d\theta_2$ where $F$ is some smooth
positive function. If~$\g$ is any metric with prescribed null cone field
$T_1\Sigma\oplus{}T_2\Sigma$, we denote by
\be
\class{\g}=\{F\cdot\g|F\in{}C^\infty(\Sigma,\bbR^+_*)\}
\label{classg1}
\ee
the class of metrics  conformally equivalent to $\g$. Thus, a conformal
Lorentz structure \cite{Wei} on $\Sigma$ is equivalently defined by
$(\Sigma,\class{\g})$.

\begin{definition}\label{DiffConf}
{\rm
A diffeomorphism $\vfi$ of $(\Sigma,\g)$ is called
\textit{conformal}---we write $\vfi\in\Conf(\Sigma,\g)$---if 
\be
\vfi^*\!\g=f_\vfi\cdot{\g}
\quad
\mbox{for some\ }
f_\vfi\in{}C^\infty(\Sigma,\bbR^+_*).
\label{defConfDiff}
\ee 
The function $f_\vfi$ is called the \textit{conformal factor}
associated with $\vfi$.}
\end{definition}

\begin{remark}
{\rm
Definition \ref{DiffConf} is general and holds in the
Riemannian case. It is, for instance, well known that $\Conf(H^2)=\PSL2$.
In the Lorentzian case, the conformal
group of $H^{1,1}$ is, however, infinite dimensional; more
precisely, we will see that $\Conf(H^{1,1})=\Diff(\S1)$. 
}
\end{remark}

\subsection{Conformal geometry of the Lorentz hyperboloid}

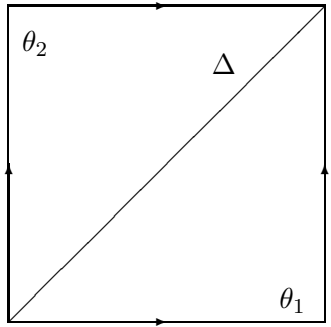
\begin{figure}[h]
\begin{center}
\begin{picture}(7,7)
\put(0,0){\framebox(7,7)}
\put(0,0){\line(1,1){7}}
\put(0,0){\vector(0,1){3.5}}
\put(0,0){\vector(1,0){3.5}}
\put(0,7){\vector(1,0){3.5}}
\put(7,0){\vector(0,1){3.5}}
\put(6,.3){$\theta_1$}
\put(.3,6){$\theta_2$}
\put(4.5,5.5){$\Delta$}
\end{picture}
\end{center}
\caption{\label{T2MinusDelta}The hyperboloid}
\end{figure}

We have seen (\ref{KostantSternberg2}) that the global intersection
properties of the rulings of the hyperboloid yield (see
Figure~\ref{T2MinusDelta})
\be
H=\bbT\x\bbT-\Delta
\label{KostantSternberg2bis}
\ee
whose metric (\ref{KostantSternberg1},\ref{theMetric}) is
given by
\be
\g_1 =
\frac{4\,d\theta_1d\theta_2}
{\left|e^{i\theta_1}-e^{i\theta_2}\right|^2}.
\label{g1}
\ee

In view of the previous definitions \ref{DefConf} and \ref{DiffConf}, any
conformal (grid-preserving) diffeomorphism $\vfi$ of a Lorentz surface
$(\Sigma,\g)$ is, locally, of the form $\vfi_1\times\vfi_2$ where $\vfi_j
\in \Diff(N_j)$. A (global) conformal diffeomorphism of~$\Sigma$ must
preserve the two foliations by lines. 

In our case, such a transformation of $\bbT^2-\Delta$ must preserve not
only the meridians and parallels of $\bbT^2$, but the diagonal $\Delta$ as
well. Therefore, $\vfi_1(\theta) = \vfi_2(\theta)$ for all
$\theta\in\bbT$, whence the
\begin{proposition}[\cite{KS2}]  
There exists a canonical isomorphism
$$
\Diff(\Delta)\stackrel{\cong}\longrightarrow\Conf(H)
$$
given by the diagonal map: $\vfi\longmapsto\vfi\times\vfi$.
\end{proposition}

\goodbreak

Let us recall the
\begin{theorem}[\cite{KS2}]\label{KS.th}
(i) Let $\vfi\in\Diff_+(\S1)\cong\Conf_+(H)$ be given. Then
$f_\vfi=(\vfi^*\!\g_1)/\g_1\longrightarrow1$ as one tends to the conformal
boundary
$\Delta$.

(ii) The conformal factor $f_\vfi$ extends smoothly to
$H\cup\Delta=\bbT^2$ and has, moreover, $\Delta$ as its critical set.

(iii) One has $\Hess(f_\vfi)|\Delta=\frac{1}{3}\tS(\vfi)$
where
\be
\tS(\vfi)
=
S(\vfi) +
\frac{1}{2}\left(\vfi'(\theta)^2-1\right) d\theta^2.
\label{SchwKS}
\ee

(iv) The Schwarzian $\tS(\vfi)$ completely determines $f_\vfi$.
\end{theorem} 
Our proof proceeds as follows. Comparison with the definition
(\ref{theMetric}) of the metric $\g_1$ on the hyperboloid $H$ in terms of
the cross-ratio prompts the following computation. Given
any $\vfi\in\Diffp(\bbT)$ viewed as a conformal diffeomorphism
(\ref{defConfDiff}) of $(H,\g_1)$, apply the Cartan formula
(\ref{SchwarzCartan}) in the case of a diffeomorphism of the circle
$\bbT$, and get
\begin{eqnarray}
\frac{(\vfi^*\!\g_1)(\theta_1,\theta_2)}{\g_1(\theta_1,\theta_2)}-1
&=&
\label{fminusone1}
f_\vfi(\theta_1,\theta_2)-1\\
&=&
\frac{1}{6}
S(\tilde\vfi)(e^{i\theta})
\left(e^{i\theta_1}-e^{i\theta_2}\right)^2+\cdots
\label{SchwarzHat1}
\end{eqnarray}
where $\tilde\vfi(e^{i\theta})=e^{i\vfi(\theta)}$ and
$\theta_j\longrightarrow\theta$ for $j=1,2$.
A tedious calculation using (\ref{SchwarzCartan}) leads to
\begin{lemma}
If $\tilde\vfi\in\Diffp(\bbT)$ is represented by\footnote{
We denote by $\Diff_{2\pi\bbZ}(\bbR)$ the universal covering of
$\Diffp(\S1)$, i.e., the group of those diffeomorphisms $\vfi$ of
$\bbR$ such that $\vfi(\theta+2\pi)=\vfi(\theta)+2\pi$.
}  
$\vfi\in\Diff_{2\pi\bbZ}(\bbR)$, one has
\be
S(\tilde\vfi)(e^{i\theta})
=
-\left(
S(\vfi)(\theta)
+
\frac{1}{2}\left(\vfi'(\theta)^2-1\right)
\right)
e^{-2i\theta}.
\label{SchwarzHat}
\ee
\end{lemma}
From (\ref{fminusone1})--(\ref{SchwarzHat}) one obtains
\be
f_\vfi(\theta_1,\theta_2)-1
=
\frac{1}{6}
\left(
S(\vfi)(\theta)
+
\frac{1}{2}\left(\vfi'(\theta)^2-1\right)
\right)
(\theta_1-\theta_2)^2
+\cdots
\label{SchwarzHat2}
\ee
that is, theorem 1 in \cite{KS2}. In particular, the conformal factor
$f_\vfi$ extends to the diagonal $\Delta\subset\bbT^2$ (its critical set)
and $f_\vfi|\Delta=1$, its transverse Hessian being related to the
modified Schwarzian derivative (see (\ref{SchwarzHat2})) by
$
\Hess(f_\vfi)
=
\frac{1}{3}\tS(\vfi)
$.
The fourth item of theorem \ref{KS.th} will be a consequence of Theorem
\ref{main}.

\goodbreak

\medskip

We are thus led to the
\begin{theorem}\label{FirstTheorem}
(i) Given any $\vfi\in\Conf_+(H)$ of $H=\bbT^2-\Delta$ and
$c\neq0$, the twice-symmetric tensor field $\vfi^*\!\g_c-\g_c$ of $H$
extends to null infinity~$\Delta$ and defines a non trivial
$1$-cocycle
\be
S_c:\vfi
\longmapsto 
\frac{3}{2}
\left(\strut\vfi^*\!\g_c - \g_c\right)\!|\Delta
\label{TheCocycle}
\ee
of $\Diffp(\bbT)$ with values in the module $\cQ(\bbT)$ of quadratic
dif\-ferentials of the circle, given by the (modified) Schwarzian
derivative (\ref{SchwKS}):
\be
S_c=c\,\tS.
\label{Schwarzc}
\ee

(ii) There holds
$H^1(\Diffp(\bbT),\cQ(\bbT))=\bbR\left[S_1\right]$.
\end{theorem}
\textit{Proof:} 
From the formul\ae\ (\ref{SchwarzHat2}) and (\ref{g1}) one immediately
gets
\begin{eqnarray*}
\left(\vfi^*\!\g_1 - \g_1\right)|\Delta
&=&
\left({\textstyle\frac{2}{3}}\,\tS(\vfi)(\theta)(\theta_1-\theta_2)^2
d\theta_1d\theta_2|e^{i\theta_1}-e^{i\theta_2}|^{-2}+\cdots\right)
\!\Big|\Delta\\[6pt]
&=&
\left({\textstyle\frac{2}{3}}\,\tS(\vfi)(\theta)d\theta_1d\theta_2+\cdots\right)
\!\Big|\Delta\\[6pt]
&=&
\frac{2}{3}\,\tS(\vfi)(\theta)\,d\theta^2\\[6pt]
&=&
\frac{2}{3}\,\tS(\vfi).
\end{eqnarray*}
Then (\ref{Schwarzc}) is clear by (\ref{KostantSternberg1}) and
(\ref{cneq0}).

At last, part (ii) follows immediately from the knowledge that
$H^1(\Diffp(\bbT),\cQ(\bbT))$ is $1$-dimensional \cite{GF,Fuchs} and
generated by the class of the Schwarzian.
\QED

\medskip
\begin{remark}
\rm{
The cocycle $\vfi\longmapsto\vfi^*\!\g_c - \g_c$ of $\Conf_+(H)$ with
values in the space of twice-covariant symmetric tensor fields is obviously
trivial. Non triviality of the cocycle (\ref{TheCocycle}) quite remarkably
stems from the ``restriction'' of the latter to null infinity
$\Delta$.\footnote{This observation is due to Valentin Ovsienko.}  
}
\end{remark}

\begin{proposition}
The group of direct isometries of the hyperboloid is
\be
\Isom_+(H,\g_c)=\ker(S_c)\cong\PSL2.
\label{IsomCurved}
\ee
\end{proposition}
\textit{Proof:} 
Using (\ref{TheCocycle}), we find that the group
$\Isom_+(H,\g_c)\subset\Diffp(\S1)$ of direct isometries is clearly a
subgroup of $\ker({S_c})\cong\PSL2$. Conversely, for any
$\vfi\in\ker({S_c})$, and thanks to (\ref{SchwarzHat1}), the conformal
factor in (\ref{fminusone1}) is $f_\vfi=1$, i.e.,
$\vfi\in\Isom_+(H,\g_c)$.
\QED

\medskip
Theorem \ref{FirstTheorem} still holds true for the $\PSL2$-invariant
metric (\ref{generalg1}) on $\cS\x\cS-\Delta$. In fact, a calculation akin to
that of (\ref{fminusone1},\ref{SchwarzHat1}) leads to
\begin{proposition}\label{TheGeneralCocycleTh}
Given any $\vfi\in\Conf_+(\cH)$ of $\cH=\cS\x\cS-\Delta$ where $\cS$ is
endowed with the projective structure (\ref{T}) or (\ref{P1}), one has
\be
\bS(\vfi)=
S_1(\vfi)
=
\frac{3}{2}
\left(\strut\vfi^*\!\g_1 - \g_1\right)\!|\Delta
\label{TheGeneralCocycle}
\ee
where the metric $\g_1$ on $\cH$ is given by (\ref{generalg1}) and the
universal Schwarzian~$\bS$ by
(\ref{generalSchwarzian},\ref{generalSchwarzianLoc}).
\end{proposition}

\goodbreak

\subsection{Conformal geometry of the flat cylinder}

Let us envisage, for a moment, the flat induced
Lorentz metric
\be
\g_0=d\theta_1d\theta_2
\label{g0}
\ee
on the cylinder $H=\bbT^2-\Delta$. (A non significant constant factor might
be introduced in the definition (\ref{g0}) of $\g_0$.)

In this special case, the $\Diffp(\S1)$-cocycle $S_0$ defined, in the
same manner as in (\ref{TheCocycle}), by
\be
S_0(\vfi)=\left(\strut\vfi^*\!\g_0-\g_0\right)\vert\Delta
\label{S0}
\ee
is, plainly, a coboundary since $\g_0$ admits a prolongation to $\Delta$.
We, indeed, have
$S_0(\vfi)(\theta)=(\vfi'(\theta)^2-1)\,d\theta^2$. Notice that flatness
of the metric is now related to triviality of the associated cocycle.

\begin{proposition}
The group of direct isometries of the flat cylinder is
\be
\Isom_+(H,\g_0)=\ker(S_0)\cong\bbT.
\label{isomzero}
\ee
\end{proposition}
\textit{Proof:} 
Solving $\vfi^*\!\g_0=\g_0$ and $\vfi'(\theta)>0$ gives 
$\vfi(\theta)=\theta+t$ with $t\in\bbT$, that is $\vfi\in\ker(S_0)$. 
\QED

\section{Symplectic structure on conformal classes of metrics on
$\cS\x\cS-\Delta$}
\label{SymplecticStructuresOnConformalClassesOfMetrics}

We analyze, in this section, the structure of the conformal classes of the
previously introduced metrics $\g_c$ and $\g_0$ on the ``hyperboloid''
$\cH$ and relate them to the generic coadjoint orbits \cite{Kiri} in
the regular dual of the Virasoro group. It should be recalled that the
conformal class of $\g_1$ has first been identified with the homogeneous
space $\Diffp(\S1)/\PSL2$ in \cite{KS1}.

\subsection{Homogeneous space $\Diffp(\cS)/\PSL2$}
\label{DiffOverPSL2}

\subsubsection{Conformal classes of curved metrics}

Consider first the curved case. If $c\neq0$, denote by $M_c$ the
space of metrics on $\cH=\cS\x\cS-\Delta$ related to $\g_c=c\,\g_1$ 
(\ref{KostantSternberg1}) by a conformal diffeomorphism (see
(\ref{classg1})), \textit{viz.}
$$
M_c=\{\g\in\class{\g_1}\,|\,\g=\vfi^*\!\g_c,\vfi\in\Conf_+(\cH)\}.
$$


\begin{figure}[h]
\begin{center}
\includegraphics[scale=0.53]{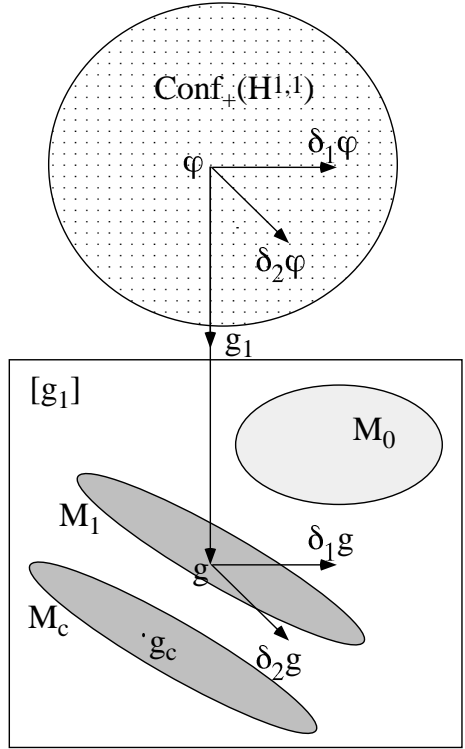}
\caption{\label{ConfClassesEPSF}
\textit{The conformal classes of metrics on $\bbT^2-\Delta$}}
\end{center}
\end{figure}


These classes $M_c$ of metrics (see Figure \ref{ConfClassesEPSF}) turn out
to have a symplectic structure of their own.
\begin{theorem}\label{ConfOrbitsTh}
If $c\neq0$, the homogeneous space
\begin{eqnarray*}
M_c
&=&
\Im(\vfi\longmapsto\vfi^*\!\g_c)\\
&\cong&
\Conf_+(\cH)/\Isom_+(\cH,\g_c)
\end{eqnarray*}
is endowed with (weak) symplectic structure $\omega_c$
which reads
\be
\omega_c(\delta_1\!\g,\delta_2\!\g)
=
\frac{3}{2}
\int_\Delta{\strut\!i_{\xi_1}L_{\xi_2}\!\g}
\label{theSymplecticStructure}
\ee
where $\delta_j\!\g=L_{\xi_j}\!\g$ with $\xi_j\in\Vect(\cS)$.
\end{theorem}
\textit{Proof:} 
From (\ref{omega1}) below, $\omega_c$ is, indeed, skew-symmetric in its
arguments. It is, clearly, also closed. We then have
$\delta_2\!\g_c\in\ker\omega_c$ iff
$\omega_c(\delta_1\!\g_c,\delta_2\!\g_c)=0$ for all $\xi_1\in\Vect(\cS)$,
i.e., iff $L_{\xi_2}\!\g_c\!|\Delta=0$, that is iff $\delta_2\!\g_c=0$
in view of (\ref{TheCocycle}) and (\ref{IsomCurved}).
\QED

\medskip
We will prove that $M_c$ is symplectomorphic to a
Kirillov-Segal-Witten $\Diff_+(\cS)$-orbit \cite{Kiri,Segal,Wit} for the
affine coadjoint (anti-)action $\Coad_\Theta$ on $\cQ(\cS)$ defined by
\be
\Coad_\Theta(\vfi)q = \Coad(\vfi)q+\Theta(\vfi)
\label{CoadS}
\ee
where the $\Diff_+(\cS)$-coadjoint (anti-)action reads
\be
\Coad(\vfi)q=\vfi^*q
\label{DefCoad}
\ee
and where $\Theta$, a $1$-cocycle of $\Diff_+(\cS)$ with values in
$\cQ(\cS)$, is a particular Souriau cocycle \cite{JMS}.

\goodbreak

\subsubsection{Intermezzo}

This technical section presents the standard $\Diffp(\cS)$-cocycles
in a guise adapted to any projective structure (\ref{T},\ref{P1}) on the
circle~$\cS$. 

Consider the line element 
$$
\lambda=\Phi^*d\theta
$$ 
on $\cS$ associated with the developing map
$\Phi\in\Diff_\mathrm{loc}(\bbR,\RP1)$. Actually, $\lambda$ is a
$\pi_1(\cS)$-invariant line-element of $\bbR$ which therefore descends to
$\cS$.

Let $\vfi$ be a representative in
$\Diff_{\pi_1(\cS)}(\bbR)$ of a diffeomorphism of $\cS$ and let
$\tilde\vfi=\Phi\circ\vfi\circ\Phi^{-1}$ denote the diffeomorphism it
induces on $\RP1$.

\begin{proposition}
(i) The \textrm{Euclidean cocycle} $\bE(\vfi)=\Phi^*E(\tilde\vfi)$ where
$E(\tilde\vfi)=\log((\tilde\vfi^*d\theta)/d\theta)$ reads
\be
\bE(\vfi)=\log\left(\frac{\vfi^*\lambda}{\lambda}\right).
\label{Euclide}
\ee

\goodbreak

(ii) The \textrm{affine cocycle} $\bA(\vfi)=\Phi^*dE(\tilde\vfi)$ is then
\be
\bA(\vfi)=d\bE(\vfi).
\label{Affine}
\ee

(iii) The \textrm{Schwarzian cocycle}
$\bS(\vfi)=\Phi^*S(\tilde\vfi)$ (see
(\ref{generalSchwarzian},\ref{generalSchwarzianLoc})) retains the form
\be
\bS(\vfi)
=
\lambda\,d\left(\frac{\bA(\vfi)}{\lambda}\right)
-\frac{1}{2}\bA(\vfi)^2.
\label{NiceSchwarzian}
\ee
\end{proposition}
\textit{Proof:} 
We easily prove (iii) by noticing that the Schwarzian
(\ref{TheFormule}) can be written in term of the affine coordinate
$\theta$ of $\RP1$ as
$$
S(\tilde\vfi)
=
d\theta\,d\left(\frac{\tilde\vfi''(\theta)}{\tilde\vfi'(\theta)}\right)
-
\frac{1}{2}\left(
\frac{\tilde\vfi''(\theta)}{\tilde\vfi'(\theta)}\,d\theta
\right)^2
$$ 
and the affine cocycle
as $A(\tilde\vfi)=(\tilde\vfi''(\theta)/\tilde\vfi'(\theta))\,d\theta$.
\QED

\medskip
For example, the $\Diffp(\S1)$-Schwarzian in angular
coordinate is recovered with
$\Phi$ as in (\ref{T}); one finds
$$
\bS(\vfi)(\theta)=\tS(\vfi)(\theta),
$$
i.e., the modified Schwarzian derivative (\ref{SchwKS}). See also
\cite{Segal}.

\begin{proposition}\label{TheInfinitesimalCocycle}
The infinitesimal Schwarzian takes either forms
$$
\bs(\xi)=s_1(\xi)
$$
for any $\xi\in\Vect(\cS)$
with\footnote{Recall that $\Div\xi=(L_\xi\lambda)/\lambda$.}
\be
\bs(\xi)
=
\lambda\,d\left(\frac{d\Div\xi}{\lambda}\right)
\label{gf}
\ee
and
$$
s_1(\xi)
=
\threehalves{}\left(L_\xi\g_1\strut\right)\!|\Delta.
$$
\end{proposition}
\textit{Proof:} 
This follows clearly from (\ref{TheGeneralCocycle}) and
(\ref{NiceSchwarzian}). 
\QED
\begin{remark}
{\rm 
In local affine coordinate on $\RP1$, the infinitesimal Schwarzian
(\ref{gf}) of $\xi=\xi(t)\partial/\partial{t}$ retains the familiar form
$$
\bs(\xi)=\xi'''(t)\,dt^2.
$$
}
\end{remark}

\subsubsection{A Virasoro orbit}

With these preparations, let us formulate the
\begin{proposition}\label{dalphaProp}
Endow $\Diff_+(\cS)$ with the $1$-form $\alpha$
defined
by
\be
\alpha(\delta\vfi)
=
\frac{1}{2}\int_{\cS}{\!\bA(\vfi)\delta{}\bE(\vfi)}
\label{alpha}
\ee
where 
$\delta\vfi=\delta(\vfi\circ\psi)$ with
$\delta\psi=\xi\in\Vect(\cS)$ at $\psi=\id$.

(i) The exterior derivative of $\alpha$ is given, for
$\xi_1,\xi_2\in\Vect(\cS)$, by
\be
d\alpha(\delta_1\vfi,\delta_2\vfi)
=
\int_{\cS}{\!\bS(\vfi)([\xi_1,\xi_2])}
+
\underbrace{
\int_{\cS}{\!d(\Div\xi_1)\,\Div\xi_2}.
}_{\GF(\xi_1,\xi_2)}
\label{dalpha}
\ee

(ii) If $\sigma$ denotes the canonical symplectic structure of the
$\Diffp(\cS)$-affine coadjoint orbit $\cO$ of the origin with Souriau
cocycle $\bS$ (see (\ref{CoadS})), namely if
\begin{eqnarray}
\cO &=& \Im(\bS)\\
&\cong&\Diffp(\cS))/\PSL2
\label{O}
\end{eqnarray}
then
\be
d\alpha=\bS^*\sigma.
\label{KKS}
\ee
\end{proposition}
\textit{Proof:} 
Since
$
d\alpha(\delta_1\vfi,\delta_2\vfi)
=
\half\int_{\cS}{d(\delta_1\bE(\vfi))\delta_2\bE(\vfi)}
-
\half\int_{\cS}{d(\delta_2\bE(\vfi))\delta_1\bE(\vfi)}
$
let us first remark that 
$$
\delta_j\bE(\vfi)=\bA(\vfi)(\xi_j)+\Div\xi_j
$$
with the above notation. If we posit for convenience $a=\bA/\lambda$, and
note that
$\lambda(\xi_1)\Div\xi_2-\lambda(\xi_2)\Div\xi_1=\lambda([\xi_1,\xi_2])$,
a lengthy calculation then leads to
$$
d\alpha(\delta_1\vfi,\delta_2\vfi)
=
\int_{\cS}{(da-\half{}a^2\lambda)\lambda([\xi_1,\xi_2])}
+
\int_{\cS}{d(\Div\xi_1)\Div\xi_2}.
$$
Whence the sought equation (\ref{dalpha}).

Now, the affine coadjoint orbit of $q_1\in\cQ(\cS)$ given
by the action~(\ref{CoadS}) carries a canonical
symplectic structure $\sigma$ which reads \cite{JMS}:
\be
\sigma(\delta_1q,\delta_2q)
=
\langle{}q,[\xi_1,\xi_2]\rangle + f(\xi_1,\xi_2)
\label{Souriau}
\ee
at $q=\Coad_\Theta(\vfi)q_1$; here $f\in{}Z^2(\Vect(\cS),\bbR)$ is the
derivative of the group-cocycle
$\Theta\in{}Z^1(\Diffp(\cS),\cQ(\cS))$ at the identity. The
expression~(\ref{dalpha}) of $d\alpha$ clearly matches that of
$\sigma$ (\ref{Souriau}) with $q_1=0$, $\Theta=\bS$ and $f=\GF$
where the Gelfand-Fuchs cocycle \cite{GF} reads
\be
\GF(\xi_1,\xi_2)
=
-\int_\cS{\!\bs(\xi_1)(\xi_2)}
\label{GelfandFuchs}
\ee
according to (\ref{gf}).
\QED

\goodbreak

Our main result is then given by
\begin{theorem}\label{main}
The map
\be
J_c:
g\longmapsto
\frac{3}{2}\left(\strut\g-\g_c\right)\!\big|\Delta
\label{Jc}
\ee
establishes a symplectomorphism\footnote{It is the momentum
map of the hamiltonian action of $\Conf_+(\cH)$ on
$(M_c,\omega_c)$.}
\be
J_c:(M_c,\omega_c)\longrightarrow(\cO_c,\sigma_c)
\label{ourSymplecto}
\ee
between the metrics of $\cH=\cS\x\cS-\Delta$ conformally
related to $\g_c$ and the affine coadjoint orbit $\cO_c=c\cdot\cO$
(see (\ref{O})) with central charge $c$, the inverse curvature
(\ref{curvature}).
\end{theorem} 
\textit{Proof:} 
Let us denote by $\g_c:\Conf_+(\cH)\longrightarrow{}M_c$ the orbital map
and let us put $\g=\g_c(\vfi)=\vfi^*\!\g_c$. We find, using
(\ref{theSymplecticStructure}),
\begin{eqnarray*}
\omega_1(\delta_1\!\g,\delta_2\!\g)
&=&
\threehalves
\int_\Delta{\!i_{\xi_1}L_{\xi_2}(\g-\g_1)}
+
\threehalves
\int_\Delta{\!i_{\xi_1}L_{\xi_2}(\g_1)}\\[6pt]
&=&
\threehalves
\int_\Delta{\!(\g-\g_1)([\xi_1,\xi_2])}
+
\threehalves
\int_\Delta{\!i_{\xi_1}L_{\xi_2}(\g_1)}\\[6pt]
&=&
\int_\Delta{\!S_1(\vfi)([\xi_1,\xi_2])}
-
\int_\Delta{\!s_1(\xi_1)(\xi_2)}\\[6pt]
&=&
\int_\Delta{\!\bS(\vfi)([\xi_1,\xi_2])}
-
\int_\Delta{\!\bs(\xi_1)(\xi_2)}
\end{eqnarray*}
with the help of Propositions \ref{TheGeneralCocycleTh} and 
\ref{TheInfinitesimalCocycle}. 
Note that we have taken into account the skew-sym\-metry
of the Gelfand-Fuchs cocycle introduced in (\ref{dalpha}) and
(\ref{GelfandFuchs}). One thus gets
\be
\omega_1(\delta_1\!\g,\delta_2\!\g)
=
\langle\bS(\vfi),[\xi_1,\xi_2]\rangle+\GF(\xi_1,\xi_2)
\label{omega1}
\ee
and, since $\g_c=c\,\g_1$,
$$
\omega_c=c\,\omega_1.
$$
Thanks to (\ref{dalpha}) and (\ref{KKS}), one can claim that
\begin{eqnarray*}
d\alpha
&=&
\g_1^*\omega_1\\
&=&\bS^*\sigma.
\end{eqnarray*}
At last, this clearly entails
$$
\omega_c=J_c^*\sigma_c
$$
where $\sigma_c=c\,\sigma$ is the canonical symplectic structure on
$\cO_c$.
\QED

\goodbreak

The following diagram summarizes our claim.
$$
\diagram
{
\Conf_+(\cH) & \hfl{\cong}{} & \Diffp(\cS)\cr
\vfl{\g_c}{} && \vfl{}{c\,\bS}\cr
M_c & \hfl{\cong}{J_c} & \cO_c\cr
}
$$

\subsection{Homogeneous space $\Diffp(\cS)/\S1$}
\label{DiffOverS1}

Consider then the flat case (\ref{g0}) and introduce the space $M_0$
of metrics (see Figure \ref{ConfClassesEPSF}) on $\cH=\cS\x\cS-\Delta$ related
to $\g_0$ by a conformal diffeomorphism, \textit{viz.}
$$
M_0=\{\g\in\class{\g_1}\,|\,\g=\vfi^*\!\g_0,\vfi\in\Conf_+(\cH)\}.
$$
\begin{theorem}
The homogeneous space 
\begin{eqnarray*}
M_0
&=&
\Im(\vfi\longmapsto\vfi^*\!\g_0)\\
&\cong&
\Conf_+(\cH)/\Isom_+(\cH,\g_0)
\end{eqnarray*}
is endowed with a (weak) symplectic structure $\omega_0$ which reads
\begin{eqnarray}
\omega_0(\delta_1\!\g,\delta_2\!\g)
&=&\label{omegazero}
\int_\Delta{\strut\!i_{\xi_1}L_{\xi_2}\!\g}\\[6pt]
&=&
\int_\Delta{\strut\!\!\g([\xi_1,\xi_2])}
\label{omegazeroBis}
\end{eqnarray}
where $\delta_j\!\g=L_{\xi_j}\!\g$ with $\xi_j\in\Vect(\cS)$.
\end{theorem}
\textit{Proof:} 
Since $\g_0$ can be prolongated to $\Delta$, (\ref{omegazero}) may
be rewritten as (\ref{omegazeroBis}) which is manifestly
skew-symmetric in its arguments. The closed $2$-form $\omega_0$ is weakly
non degenerate as
$\delta_2\!\g\in\ker\omega_0$ iff $L_{\xi_2}\!\g\!|\Delta=0$, i.e.
$\delta_2\!\g=0$ in view of (\ref{S0}) and (\ref{isomzero}).
\QED

\medskip
In fact, $M_0$ is symplectomorphic to a $\Diffp(\cS)$-coadjoint
orbit \cite{Kiri} as shown below.

\goodbreak

Let us consider the following quadratic differential
\be
q_0=\g_0\!|\Delta\in\cQ(\cS)
\label{q0}
\ee
so that the $\Diffp(\cS)$-coadjoint (anti-)action\footnote{We, indeed,
have $\Coad(\vfi)(q_0)=(q_0\circ\Ad)(\vfi)$ for all $\vfi\in\Diffp(\cS)$.} 
$\Coad$ given by (see (\ref{DefCoad})) 
$\Coad(\vfi):q_0\longmapsto{}q=\vfi^*q_0$, reads according to (\ref{S0}):
\be
q=q_0+S_0(\vfi).
\label{Coad}
\ee

\begin{proposition}\label{dalphazeroProp}
Endow $\Diff_+(\cS)$ with the $1$-form $\alpha_0$ defined by
$$
\alpha_0(\delta\vfi)
=
-\int_{\cS}{\!(\vfi^*q_0)(\xi)}
$$
where, again, 
$\delta\vfi=\delta(\vfi\circ\psi)$ with
$\delta\psi=\xi\in\Vect(\cS)$ at $\psi=\id$.

(i) We have, for any $\xi_1,\xi_2\in\Vect(\cS)$,
$$
d\alpha_0(\delta_1\vfi,\delta_2\vfi)
=
\int_{\cS}{\!(\vfi^*q_0)([\xi_1,\xi_2])}.
$$

(ii)  The $\Diff_+(\cS)$-coadjoint orbit through $q_0$ (\ref{q0}) is
\begin{eqnarray}
\cO_{q_0}
&=&\Im(q_0\circ\Ad)\\
&\cong&
\Diff_+(\cS)/\S1
\label{Oq0}
\end{eqnarray}
and is endowed with the symplectic $2$-form $\sigma_0$ such that
$$
d\alpha_0=(q_0\circ\Ad)^*\sigma_0.
$$
\end{proposition}
\textit{Proof:} 
If $\delta_j\vfi$ is associated with $\xi_j\in\Vect(\cS)$ at
$\vfi\in\Diffp(\cS)$, one readily finds $\delta_jq=L_{\xi_j}q$ and
$d\alpha_0(\delta_1\vfi,\delta_2\vfi)
=
-\alpha_0([\delta_1,\delta_2]\vfi)
=
\langle{}q,[\xi_1,\xi_2]\rangle
$
which descends as the canonical symplectic $2$-form $\sigma_0$ of
$\cO_{q_0}$, namely
$$
d\alpha_0(\delta_1\vfi,\delta_2\vfi)
=
\sigma_0(\delta_1q,\delta_2q).
$$
We then simply check that $\ker(d\alpha_0)$ is $1$-dimensional and
integrated by $\ker(S_0)\cong\bbT$ (see (\ref{isomzero}) and (\ref{Coad})).
\QED

\medskip
The ``flat'' counterpart of Theorem \ref{main} is now at hand.
\begin{theorem}\label{mainzero}
The map
\be
J_0:
\g\longmapsto\g\big|\Delta
\ee
establishes a symplectomorphism\footnote{It is the momentum
map of the hamiltonian action of $\Conf_+(\cH)$ on
$(M_0,\omega_0)$.}
\be
J_0:(M_0,\omega_0)\longrightarrow(\cO_{q_0},\sigma_0)
\ee
between the metrics of $\cH=\cS\x\cS-\Delta$ conformally
related to $\g_0$ and the coadjoint orbit $\cO_{q_0}$ (see \ref{Oq0}) with
zero central charge.
\end{theorem} 
\textit{Proof:} 
Clear.
\QED

\subsection{Bott-Thurston cocycle and contactomorphisms}

It is know since the work of Kirillov \cite{Kiri} that the
$\Diffp(\cS)$-homogeneous spaces we dealt with in Sections
\ref{DiffOverPSL2} and \ref{DiffOverS1} are, in fact, genuine coadjoint
orbits of the Virasoro group, $\Vir$, i.e., the $(\bbR,+)$-central
extension \cite{TW} of $\Diffp(\cS)$ that can be recovered as follows in
our setting.

Let us emphasize that the $1$-form $\alpha$ (\ref{alpha}) on $\Diff_+(\cS)$
fails to be invariant. So, let us equip
$\Diff_+(\cS)\times\bbR$ with the following ``contact'' $1$-form~$\halpha$, \textit{viz.}
\be
\halpha(\delta\vfi,\delta{t})=\alpha(\delta\vfi)+\delta{t}.
\label{hatAlpha}
\ee
Now, the $2$-form $d\halpha$ is $\Diff_+(\cS)$-invariant and plainly
descends to $M_1$ as~$\omega_1$ (see (\ref{KKS}) and
(\ref{Jc},\ref{ourSymplecto})). We now have the
\begin{proposition}
Lifting $\Diff_+(\cS)$ into the group of auto\-morphisms of
$(\Diff_+(\cS)\times\bbR,\halpha)$ yields the Virasoro group $\Vir$
with multiplication law
\be
(\vfi_1,t_1)\cdot(\vfi_2,t_2)
=
\Big(
\vfi_1\circ\vfi_2,t_1+t_2
\underbrace{
-\frac{1}{2}\int_{\cS}{\!\bE(\vfi_1\circ\vfi_2)\bA(\vfi_2)}
}_{\BT(\vfi_1,\vfi_2)}
\Big)
\label{BT}
\ee
where $\BT$ is the Bott-Thurston cocycle \cite{Bott} of
$\Diffp(\cS)\cong\Conf_+(\cH)$.
\end{proposition}
\textit{Proof:} 
Using the cocycle relation
$\bE(\vfi\circ\psi)=\psi^*\bE(\vfi)+\bE(\psi)$---see (\ref{Euclide})---and
(\ref{Affine},\ref{alpha}), one immediately finds
\begin{eqnarray*}
\alpha(\delta(\vfi\circ\psi))
&=&
\half\int_\cS{\!\psi^*(\bA(\vfi)\delta(\bE(\vfi))}
+
\half\int_\cS{\!\bA(\psi)\delta(\bE(\vfi\circ\psi))}\\[6pt]
&=&
\alpha(\delta\vfi)
+
\delta\left[
\half\int_\cS{\!\bE(\vfi\circ\psi)\bA(\psi)}
\right]
\end{eqnarray*}
for all $\vfi,\psi\in\Diffp(\cS)$. Looking for those maps
$(\vfi,t)\longmapsto(\vfi^\star,t^\star)$ such that
$\vfi^\star=\vfi\circ\psi$ and
$\halpha(\delta\vfi^\star,\delta{t^\star})=\halpha(\delta\vfi,\delta{t})$
leads to $t^\star=t+\BT(\vfi,\psi)+\const$, hence, to the group law
(\ref{BT}).
\QED

\medskip
The triple $(\bS,\GF,\BT)$ is a special instance of a general
structure that has been coined ``trilogy of the moment'' \cite{Igl}.
\begin{remark}
{\rm
It would be interesting to have a conformal interpretation of
the contact structure $\Vir/(\ker\halpha\cap\ker{d\halpha})$ above
$(M_1,\omega_1)$.
}
\end{remark}

\goodbreak
 
\section{Conclusion and outlook}\label{ConclusionAndOutlook}


This work prompts a series of more or less ambitious questions connected
with the striking analogies between conformal geometry of Lorentz surfaces
and projective geometry of conformal infinity that we have just discussed. 
It constitutes an introduction to a more detailed paper (in preparation)
where the authors wish to tackle the following problems.
\begin{enumerate}

\item Is it possible to realize any Virasoro coadjoint
orbit\footnote{Other isotropy groups are, e.g., the finite coverings of
$\PSL2$ and $1$-parameter subgroups of the form $\S1\x\bbZ_n$; see
\cite{Guieu2}.} as a conformal class of Lorentz metrics on the cylinder?
If this is so, spell out the symplectic forms in terms of the classes of
metrics; also study the relationship between the properties of an
orbit and the dynamics of the null foliations in the associated conformal
class. There exists, in fact, a map sending the space of Virasoro
orbits---modules of projective structures on the circle---to the space of
modules of Lorentzian conformal structures on the cylinder; analyze its
properties. More conceptually, given a conic $C$ in the real
projective plane, what are the links between the space of projective
structures on $C$, the space of Lorentzian structures in the exterior of $C$
and the space of Riemannian metrics in the interior of $C$?

\item The Ghys theorem \cite{Ghys,OT} states that any diffeomorphism of the
projective line has at least four points where its Schwarzian vanishes,
i.e., four points where the contact of the graph of the diffeomorphism
with its osculating hyperbola is greater than the generic one. This result
is a Lorentzian analogue of the so-called four vertices
theorem\footnote{Any closed simple curve in the plane admits at least four
points where its Euclidean curvature is critical.} for closed curves in
the Euclidean plane. In our context, the Ghys theorem would imply the
existence, for any conformal automorphism of the hyperboloid, of some
particular points where this diffeomorphism is closer than usual to an
isometry.

\item The orbit $\Diff(\S1)/\PSL2$ embeds symplectically in the universal
Teichm\"uller space $T(1)=\QS(\S1)/\PSL2$, where $\QS(\S1)$ denotes the
group of quasi-symmetric homeomorphisms of the circle~\cite{NV}. With the
help of the quantum differential calculus of Connes, it is possible to
construct extensions of the three fundamental cocycles $\bE$, $\bA$ and
$\bS$ to the group
$\QS(\S1)$~\cite{NS}. Can one construct a ``quantum analogue'' of the
Lorentzian hyperboloid whose conformal class may be identified with~$T(1)$?
\end{enumerate}

Let us finally mention two other subjects closely connected with
our problem, namely the geometry of the Wess-Zumino-Witten model \cite{FG}
and Douglas' proof of the Plateau problem revisited by Guillemin,
Kostant and Sternberg \cite{GKS}.




\end{document}